\documentclass[conference]{IEEEtran}
\ifCLASSINFOpdf
\else
\fi

\usepackage{cite}
\usepackage{hyperref}
\usepackage{url}
\usepackage{graphicx}
\usepackage{graphics} 
\usepackage{epsfig} 
\usepackage{mathptmx} 
\usepackage{mathrsfs}
\usepackage{times} 
\usepackage{amsmath} 
\usepackage{amssymb}  
\usepackage{color}
\usepackage{float}
\usepackage{algorithm}
\usepackage{algorithmic}
\usepackage[utf8]{inputenc} 
\usepackage[T1]{fontenc}    
\usepackage{booktabs}       
\usepackage{amsfonts}       
\usepackage{nicefrac}       
\usepackage{microtype}      
\usepackage{color}
\usepackage{amsmath}
\usepackage{amsthm}
\usepackage{enumitem}
\usepackage{wrapfig}
\usepackage{lipsum}
\usepackage{multirow}
\usepackage{booktabs}
\usepackage{array}
\usepackage{subcaption}
\usepackage{dsfont}
\usepackage{makecell}

\newcolumntype{P}[1]{>{\centering\arraybackslash}p{#1}}
\newcommand{\R}{\mathbb{R}}
\newcommand{\mb}{\mathbf}

\newcommand\ChangeRT[1]{\noalign{\hrule height #1}}
\newtheorem{theorem}{Theorem}

\newtheorem{lemma}{Lemma}
\newtheorem*{lemma*}{Lemma}

\DeclareMathAlphabet{\mathpzc}{OT1}{pzc}{m}{it}
\newtheoremstyle{remarkstyle}
  {3pt}                            
  {3pt}                            
  {\normalfont}                    
  {}                               
  {\bfseries}                      
  {.}                              
  {5pt plus 1pt minus 1pt}         
  {}                               

\theoremstyle{remarkstyle}

\IEEEoverridecommandlockouts
\begin{document}
\title{Optimal Voltage Control Using Online Exponential
Barrier Method
\thanks{B. Zhang is partially supported by the NSF Grant ECCS-1942326. Work done while P. Zhang was at UW.}
}
\author{\IEEEauthorblockN{Peng Zhang}
\IEEEauthorblockA{\textit{Electrical and Computer Engineering} \\
\textit{University of Wisconsin–Madison}\\
Madison, WI 53706 \\
pzhang286@wisc.edu}
\and
\IEEEauthorblockN{Baosen Zhang}
\IEEEauthorblockA{\textit{Electrical and Computer Engineering} \\
\textit{University of Washington}\\
Seattle, WA 98195 \\
zhangbao@uw.edu}
}


\maketitle


\begin{abstract}
This paper address the optimal voltage control problem of distribution systems with high penetration of inverter-based renewable energy resources, under inaccurate model information.
We propose the online exponential barrier method that explicitly leverages the online feedback from grids to enhance the robustness to model inaccuracy and incorporates the voltage constraints to maintain the safety requirements. We provide analytical results on the optimal barrier parameter selection and sufficient conditions for the safety guarantee of converged voltages. We also establish theoretical results on the exponential convergence rate with proper step-size. 
The effectiveness of the proposed framework is validated on a 56-bus radial network, where we significantly improve the robustness against model inaccuracy compared to existing methods. 
\end{abstract}

\begin{IEEEkeywords}
Nonlinear Optimization, Barrier Method, Distribution Systems, Voltage Control
\end{IEEEkeywords}

\section{Introduction}
\label{Sec:intro}

\par Voltage control plays an important role in the distribution systems, especially given that distributed energy resources (DERs) are causing larger voltage swings~\cite{farivar2012optimal,zhang2014optimal}. Since traditional mechanical devices such as tap changing transformers and shunt capacitors typically operate at slower time scales, using power outputs of DERs themselves for fast-time-scale voltage control has attracted lots of interests. 


\par One way to approach the inverter-based voltage control is to formulate them as  optimization problems ~\cite{zhang2013local, zhang2024online, dall2016optimal}. Essentially, we obtain a version of the optimal power flow (OPF) problem, which is a very well studied topic (see, e.g.~\cite{molzahn2017survey} for references). In particular, OPF for distribution systems can be solved to optimality in many settings despite its nonconvexity~\cite{low2014convex,zhang2012geometry}. In addition, there are also several linearized models~\cite{zhu2015fast,bolognani2015existence}.


\par However, this approach requires detailed parametric and topological information about the distribution systems, which might not be available
\cite{deka2023learning, ardakanian2019identification}. 
Two-stage approaches that first learn the system model then use model-based control has been proposed to overcome the lack of information~\cite{chen2020data, feng2023stability}. 
But it is not clear that how the system should be controlled during the data gathering stage and how to deal with changes in the system. Another direction of study is to implement learning and control simultaneously in real time \cite{duan2019deep, yang2019two}. In these methods, linearized power flow models are learned based on the least-square estimator \cite{fattahi2019learning} or convex body chasing \cite{yeh2022robust} through the single sample trajectory. In \cite{yeh2022robust}, the intermediate safety constraints are guaranteed by a robust control oracle, which may lead to conservative actions. 

\par More explicit online OPF methods that use physical girds as a solver have recently been developed by \cite{dall2016optimal, gan2016online, magnusson2020distributed}.
For instance, online primal-dual method is first proposed by \cite{dall2016optimal}, and then generalized to handle model inaccuracies \cite{bernstein2019real}. The solutions achieve feasibility at convergence, but not necessarily during the intermediate iterations. Alternatively, logarithmic barrier method can enforce hard constraint satisfaction \cite{davydov2023time}. But the selection of barrier weights is heuristic, which results in the slow convergence process and optimality gap.


\par In this work, we concentrate on the inverter-based optimal and safe voltage control problem under large estimation errors of distribution networks. Based on the linearized power flow equations, we propose the online exponential barrier method that achieves good performance in both optimality and safety. We analytically derive the barrier parameters to enable the online implementation, where the iterations proceed by using the observations on the measured voltage magnitudes and power injections of controlled DERs. 
Our framework not only inherits the error correction process of online methods, but also avoids the explosion around boundaries by the nature of exponential barriers.  
With the same working scheme and data requirements, the proposed framework is compatible with online system identification techniques and strikes a tradeoff between optimality and safety when large estimations remain.



\par The organization of this paper is as follows. Section \ref{Sec:problem} introduces the model of distribution systems and formulates the voltage control problem. Section \ref{Sec:alg} illustrates the main algorithm and Section \ref{Sec:theory} provides the theoretical analysis of optimality and safety. Section \ref{Sec:sim} validates the performance through extensive numerical experiments. The conclusion is drawn in Section \ref{Sec:conclu}.




\section{MODELING AND PROBLEM FORMULATION}
\label{Sec:problem}

\subsection{Modeling of distribution systems}

\par Consider a balanced three-phase  distribution system with $n+1$ nodes operating under the tree topology. By convention, the feeder is the root of the tree and labeled as bus 0. The rest of the network is represented as ${G} = ({N}, {E})$ with ${N} =\{1,2, \ldots, n\}$ as the set of non-slack nodes and ${E} \subset {N} \times {N}$ as the set of lines. Let  $v_k = |v_k| e^{j \angle v_k} \in \mathbb{C}$ and $i_k = |i_k| e^{j \angle i_k} \in \mathbb{C}$ denote the voltage phasor and injected current phasor at bus $k$, and stack them into vectors $\mathbf{v} =  [v_1 \:v_2 \:\ldots \: v_n]^T \in \mathbb{C}^n$ and $\mathbf{i} =  [i_1\: i_2\: \ldots \: i_n]^T \in \mathbb{C}^n$.

\par The voltages and currents are related as
\begin{align}
    \left[\begin{array}{c}
i_0 \\ \mathbf{i} \end{array}\right] =  \left[\begin{array}{cc}y_{0} & \mathbf{y}_{{N}}^T \\ \mathbf{y}_{{N}} & \mathbf{Y}\end{array}\right] \left[\begin{array}{c}
v_0 \\ \mathbf{v} \end{array}\right],
\label{2-0}
\end{align}
where $i_0$ and $v_0$ represent the injected current and nodal voltage at the slack bus 0. 
The entries of admittance matrix are $\mathbf{y}_{{N}} \in \mathbb{C}^n$ and  $\mathbf{Y} \in \mathbb{C}^{n\times n}$ with positive line resistance and reactance, i.e., $Y_{ik} = G_{ik} - j B_{ik}$ with $G_{ik} >0$ and $B_{ik} >0$. The inverse of $\mathbf{Y}$ is the impedance matrix $\mathbf{Z} = \mathbf{Y}^{-1} =  \mathbf{R} +j\mathbf{X} $. From \eqref{2-0} we have $\mb{v} = v_0 \mathds{1} + \mb{Z} \mb{i}$. 

\par The complex power injections at buses are defined as $ \mb{s} = \mb{p}+ j\mb{q}, \:\mb{s} \in  \mathbb{C}^n$ and satisfy the equation $\mb{s} =\text{diag}(\mb{v}) \bar{\mb{i}}$, where $\bar{\mb{i}}$ is the element-wise complex conjugate of $\mb{i}$. We use the standard 
linearized AC power flow equations found in \cite{bolognani2015existence, guggilam2016scalable}. Under the assumption of limited voltage drop in distribution networks, we obtain the following linear equations over the voltage magnitudes $
    |\mb{v}|=|v_0| \mathds{1}+\frac{1}{|v_0|} \operatorname{Re}[\mb{Z} \bar{\mb{s}}].$
    
\par We further shift the system by the constant $|v_0|$ and concentrate on the voltage deviations. For the notational simplicity, we introduce $x_k = |v_k| - |v_0| \in \mathbb{R}$ as the normalized voltage deviation at bus $k$ and vector $\mb{x} \in \mathbb{R}^n$ for $n$ buses. Consequently, we obtain the compact form of linearized AC power flow model
\begin{align}
    \mb{x} = \mb{R} \mb{p} + \mb{X} \mb{q}. 
\end{align}
where $\mb{R}$ and $\mb{X}$ are symmetric, positive definite and positive matrices, i.e., with all positive elements \cite{farivar2013equilibrium}. 


\par We consider the PV inverters installed with the capability of providing both active and reactive power. The maximum available active power is denoted by $\mb{p_{av}} \in \mathbb{R}^n$, while we assume that no reactive power is injected or consumed without control implemented. Let $\mb{u} = (\mb{u^p}, \mb{u^q}) \in \mathbb{R}^{2n}$ as the adjusted power of inverters. Then, the net active and reactive power injections can be decomposed as
\begin{align}
    \mb{x} = \mb{R} \mb{u^p} + \mb{X} \mb{u^q} + \mb{e},
\end{align}
where $\mb{e} = \mb{R}(\mb{p_{av}} - \mb{p_{e}}) + \mb{X}(-\mb{q_{e}}) \in \R^n$ represents the voltage drop caused by the original solar energy and exogenous loads $\mb{p}_e$ and $\mb{q}_e$.

\subsection{Optimal voltage control}

In this paper, the optimal voltage control problem is to maintain nodal voltages within the predefined safety intervals, typically as $[\mb{\underline{x}}, \mb{\bar{x}}] \subset \mathbb{R}^n$, by solving for the inverter actions that minimize control efforts and solar energy curtailments. We consider the physical capacity limits over all inverters, represented by $\mb{u} \in [\mb{\underline{u}}, \mb{\bar{u}}] \subset \mathbb{R}^{2n}$. Formally, we are interested in the following optimization problem
\begin{subequations}
\begin{align}
     \underset{\mb{u} \in \mathbb{R}^{2n}}{\text{minimize}}
    \quad & c(\mb{u}) \\
    \text{subject to} \quad  & \mb{\underline{x}} \leq \mb{x} \leq \mb{\bar{x}} \label{op}\\
    & \mb{x} = \mb{R} \mb{u^p} + \mb{X} \mb{u^q} + \mb{e} \\
    & \mb{\underline{u}} \leq \mb{u} \leq \mb{\bar{u}}.
    \label{2-2}
\end{align}
\end{subequations}
We further assume that the exogenous loads, $\mb{p_{e}}$ and $\mb{q_{e}}$, are well controlled by other techniques and operate within the operational constraints in \eqref{op}, such that \eqref{2-2} is always feasible. Then, we can narrow our focus on the overvoltage problem aroused by the  high penetration of solar panels. We reformulate the  addressed optimization problem as follows
\begin{subequations}
\begin{align}
     \underset{\mb{u} \in \mathbb{R}^{2n}}{\text{minimize}}
    \quad & c(\mb{u}) \label{2-3}\\
    \text{subject to} \quad  & \mb{Bu} + \mb{e} \leq \mb{\bar{x}} \label{2-4}\\
    & \mb{\underline{u}} \leq \mb{u} \leq \mb{\bar{u}},\label{2-5}
\end{align}
\end{subequations}
where $\mb{B} = [\mb{R\quad X}] \in \R^{n \times 2n}$.

\par Under the quadratic cost function, e.g. $c(\mb{u}) = \frac{1}{2}\mb{u^T Q u}$ for some $\mb{Q} \succeq 0$, the above optimization is a Linearly Constrained Quadratic Problem (LCQP) and can be efficiently solved if $\mb{B}$ is known. The goal of this paper is to solve it when $\mb{B}$ is not exactly known by leveraging real-time observations from physical grids. 

\section{Algorithm}
\label{Sec:alg}
This section proposes the online exponential barrier method for the optimal and safe voltage control under inaccurate models of distribution systems. In the closed-loop setting, the optimization-based controller updates the active and reactive power injections of inverters based on the online feedback from physical measurements.

\par We summarize the detailed steps in Algorithm \ref{alg}. In particular, to maintain the safety requirements on voltages, we integrate constraints \eqref{2-4} into the objective function as exponential barriers and obtain the following augmented cost
\begin{align}
    \tilde{c}(\mb{u}| \hat{\alpha}, \hat{\mb{B}}) = c(\mb{u}) + \sum_i^n \frac{\hat{\alpha}_i}{\beta_i} e^{\beta_i(\mb{\hat{b}}_i^T \mb{u} + \mb{e}_i - \mb{\bar{x}}_i)},
    \label{aug_cost}
\end{align}
where $\mb{\hat{b}}_i^T$ is $i$-th row of the imprecise estimation $\mb{\hat{B}}$, and $\hat{\alpha}_i$, $\beta_i, \: i \in N$ are corresponding barrier weights. Then, to deal with model inaccuracies, we adjust inverter setpoints based on the gradient calculated by the latest voltage measurements  $\mb{x}_i = \mb{b}_i^T \mb{u} + \mb{e}_i, \: i \in N$, as defined in the following vector, 
\begin{align}
    F(\mb{u} | \hat{\alpha}) := \nabla_{\mb{u}} {c}(\mb{u}) + \sum_{i}^n \hat{\alpha}_i e^{\beta_i (\mb{x}_i - \mb{\bar{x}}_i)} \hat{\mb{b}}_i.
    \label{gradient}
\end{align}
\par The key insights lie in the adoption of exponential barriers. Compared to the online primal-dual method in \cite{dall2016optimal}, our framework drives voltages within the safety limits through the significantly increasing gradient near the boundaries. In addition, exponential terms enable the analysis of voltage dynamics on the boundaries, instead of suffering from the explosion under logarithmic barriers. This property also allows us to directly compute the barrier parameters (see next section) and circumvents the double-loop process of logarithmic barrier method which relies on shrinking weights.

\begin{algorithm}
\caption{Online Exponential Barrier Method for Optimal Voltage Control}
\begin{algorithmic}
\STATE \textbf{Input 1}: Safety limit $\mb{\bar{x}}$, Cost coefficient $\mb{Q}$, Barrier curvature $\mb{\beta}$, Solar efficiency bound $\kappa$, Stopping criterion $K$
\STATE \textbf{Input 2}: Estimated dynamics $\mb{\hat{B}}$, Voltage deviation $\mb{x}$ under maximum solar energy $\mb{p_{av}}$, step-size $\eta$ 
    \IF{$\max_{i \in N} \mb{x}_i \geq \mb{\bar{x}}_i $}    
    \STATE \textbf{(s1)} Initialize the safe action $\mb{u}(k) = [(\kappa - 1) \mb{p_{av}} \quad \mb{0}]^T $ and observe new voltage deviation $\mb{x}(k)$ with $k=0$

    \STATE \textbf{(s2)} Compute barrier weights $\mb{\hat{\alpha}^s}$ for augmented cost $\tilde{c}$ 
    
    \FOR{$k = 0$ to $K$}
    \STATE \textbf{(s3)} Compute and implement adjusted power setpoint $\mb{u}(k+1) = \mb{u}(k) - \eta F(\mb{u}(k) | \alpha^s)$
    \STATE \textbf{(s4)} Measure voltage $\mb{x}(k+1)$ and control input $\mb{u}(k+1)$
    
    \IF{$\min_{i \in N} \mb{u}_i(k+1)= \mb{\underline{u}}_i$ or $\max_i \mb{u}_i(k+1)= \mb{\bar{u}}_i$}
    \STATE \textbf{(s5)} Activate saturation by updating weights $\mb{\hat{\alpha}^{s}}$
    \ENDIF
    
    \IF{$\arg \max_{i \in N} \mb{x}_i (k+1) \neq \arg \max_{i \in N} \mb{x}_i (k)$ }
    \STATE \textbf{(s6)} Switch attention node by updating weights $\mb{\hat{\alpha}^{s}}$
    \ENDIF
    
    \ENDFOR
    \ENDIF
\end{algorithmic}
\label{alg}
\end{algorithm}

\section{Optimality analysis and safety guarantee}
\label{Sec:theory}

This section provides the closed-form principles of barrier parameter selection. We first prove the optimality guarantee in error-free cases, and then extend to the safety condition under inaccurate models. At last, we show the exponential convergence rate and proper choice of step-size.   We make the following assumption about the cost function:
\par \noindent \textbf{Assumption 1}: The cost function $c(\mb{u})$ is quadratic with positive definite coefficients, i.e., $c(\mb{u}) = \frac{1}{2}\mb{u^TQu}$, with $\mb{Q} \succ 0 $.
\begin{theorem}
    (Optimality condition) Under Assumption 1 and model $\mb{B}$, there exists the unique barrier parameter $\mb{\alpha} \in \mathbb{R}_+^n$ with fixed $\mb{\beta} \in \mathbb{R}_+^{n}$, such that the augmented cost function $\tilde{c}(\mb{u} | \alpha, \mb{B})$ is strongly convex and its global minimum $\mb{\mb{u^*}}$ satisfies the optimality condition of Problem \eqref{2-3}-\eqref{2-4}, where $\alpha$ is chosen as the solution to 
    \begin{align}
        \begin{bmatrix}
    \mb{\mb{Q}} & \mb{B}_{{A}}^T \\
    \mb{B}_{{A}} & \mb{0}
    \end{bmatrix}
    \begin{bmatrix}
    \mb{u} \\
    \alpha_{{A}}
    \end{bmatrix}
    =
    \begin{bmatrix}
    \mb{0}\\
    \mb{\bar{x}}_{{A}} - \mb{e}_{{A}},
    \end{bmatrix},
    \label{optimal condition}
    \end{align}
    and $\alpha_j = 0$, $j = N\setminus {A}$, with ${A}$ as the set of active constraints at $\mb{u^*}$.
\end{theorem}

\begin{proof}
    For the augmented cost  $\tilde{c}(\mb{u} | \alpha, \mb{B})$ in \eqref{aug_cost}, by Lemmas \ref{uniqueness} and \ref{positive}, the linear equation \eqref{optimal condition} has the unique nonnegative solutions in $\alpha_{A}$. We have that  $\tilde{c}(\mb{u} | \alpha, \mb{B})$ is strongly convex.
    \par We assume the global minimum of  $\tilde{c}(\mb{u} | \alpha, \mb{B})$ as $\mb{u^*}$ and calculate the gradient as below
    \begin{align}
        \nabla_{\mb{u}} \tilde{c}(\mb{u} | {\alpha}, \mb{B}) = \mb{Qu} +\sum_{i=1}^n\alpha_ie^{\beta_i (\mb{b}_i^T \mb{u} + e_i - \bar{x}_i)}\mb{b}_i.
    \end{align}
    \par Under the strongly convexity of $\tilde{c}(\mb{u} | \alpha, \mb{B})$, we can follow the choice of $\mb{u}$ in \eqref{optimal condition}. At $\mb{u^*}$, we have $e^{\beta_i (\mb{b}_i^T\mb{u^*} + \mb{e}_i - \mb{\bar{x}}_i)} = e^0 = 1$ on the active constraints, and $\alpha_j = 0$ on the inactive constraints. Thus, we get that
    \begin{align}
        \nabla_{\mb{u}} c(\mb{u^*} ) + \sum_{i \in A} \alpha_i \mb{b}_i = 0, 
    \end{align}
    with $\alpha_i \geq 0, \: i \in {A}$.
    \par According the Lagrange Multiplier Theorem for linear constraints (\textit{Proposition 3.4.1} in \cite{bertsekas1997nonlinear}), we have the satisfaction of optimality condition for Problem  \eqref{2-3}-\eqref{2-4}.
\end{proof}

\begin{lemma}
    The linear equation \eqref{optimal condition} has  unique solutions for different length of active constraint set, i.e., $|{A}| \in [1, n]$.
\label{uniqueness}
\end{lemma}

\begin{proof}
We show that the matrix on the left hand side of~\eqref{optimal condition} is invertible. First consider the case of $|{A}| = 1$. 
Based on {Schur complement} and $\mb{Q} \succ 0$, the matrix is invertible.

Now suppose $|{A}| \in (1, n]$, taking the Schur complement gives $\mb{S} = - \mb{B}_{{A}} \mb{\mb{Q}}^{-1} \mb{B}_{{A}}^T= -\mb{N}^T \mb{N},$
where $\mb{N} = \mb{\mb{Q}}^{-1/2} \mb{B}_{{A}}^T \in \mathbb{R}^{2n \times |{A}|}$ has linearly independent rows. Since $|{A}| < 2n$, $-\mb{S}$ is positive definite and thus invertible. 
\end{proof}


\begin{lemma}
    The solutions of $\alpha_{{A}}$ to the linear system in \eqref{optimal condition} are nonnegative, i.e., $\alpha_i \geq 0, \; i \in {A}$.    
    \label{positive}
\end{lemma}

\begin{proof}
    Based on \eqref{optimal condition}, we first introduce the auxiliary vector $\mb{p = -\mb{Q}}\mb{u} = \mb{B}_{{A}}^T \alpha_{{A}}$. Since $\mb{\mb{Q}} \succ 0 $, we have $\mb{p}^T \mb{u}  = -\mb{u}^{T}\mb{\mb{Q}}\mb{u}\leq 0$ for all vectors $\mb{u} \in \mathbb{R}^{2n}$. We also have $\mb{b}_i^T \mb{u} = \bar{x}_i - e_i < 0, \; i \in {A}$. 
    Then, using Farkas's Lemma (\textit{Proposition B.16(d) in \cite{bertsekas1997nonlinear}}), we know that $-\mb{p}$ locates in the cone generated by vectors $\mb{b}_i, \: i \in {A}$. Thus, there exist nonnegative scalars $\mu_i$, such that $\mb{p} = \sum_{i \in {A}} \mu_i \mb{b}_i$.
    Recall that $\mb{p} = \sum_{i \in {A}} \alpha_i \mb{b}_i$, we obtain $\alpha_i = \mu_i \geq 0$, $i \in {A}$.
\end{proof}

\par The optimal barrier parameter $\alpha$ depends on the active constraints $A$. For the overvoltage problem, we can choose the node with the largest voltage deviation and continuously adjust the attention node based on the voltage measurements, as shown in Algorithm \ref{alg}. 

\par Faced with inaccurate model $\mb{\hat{B}}$, the above settings enable us to analyze dynamics around maximum voltage deviations, where we propose the theorem for barrier parameter $\alpha^s$ with guarantees that the converged point of Algorithm \ref{alg} always locates within the safety limits. Since inaccuracies affect both barrier calculation and action updates, we follow a two-step proof, where the existence of proper barrier weights is first demonstrated, then the tightened form is further adopted to guarantee the constraint satisfaction.

\par Since the physical capacities of inverters in \eqref{2-5} are considered, the variational inequality (VI) is introduced to study voltages at the convergence. We also make the following assumption.
\par \noindent \textbf{Assumption 2}: The cost function $c(\mb{u})$ is  homogeneous quadratic with diagonal matrix $\mb{\mb{Q}}$.

\par \noindent \textbf{Definition 1:} Given $[\mb{\underline{u}}, \mb{\bar{u}}] \subset \mathbb{R}^{2n}$,  a vector $\mb{u}^* \in \mathbb{R}^{2n}$ is called the VI solution to $F(\mb{u} | \alpha)$ if it  satisfies
\begin{align}
    F(\mb{u}^* | \alpha)^T (\mb{u} - \mb{u}^*) \ge 0, \quad \forall \mb{u} \in [\underline{\mb{u}}, \bar{\mb{u}}].
    \label{converged point1}
\end{align}

\begin{theorem}
    (Existence of $\alpha^\dagger$) Under Assumptions 1-2, given the set of unsaturated actions $A_u$ and attention node $i$, let $\alpha_i^\dagger$ be the solution to 
    \begin{align}
    \begin{bmatrix}
    \mb{\mb{Q}}_{A_u} & \mb{\hat{b}}_{{i,A_u}} \\
    \mb{b}_{i,A_u}^T & \mb{0}
    \end{bmatrix}
    \begin{bmatrix}
    \mb{u}_{A_u} \\
    \alpha_{i}
    \end{bmatrix}
    =
    \begin{bmatrix}
    \mb{0}\\
    \mb{\bar{x}}_{i} - \mb{e}_{i,A_u}
    \end{bmatrix},
    \label{def1}
    \end{align}
with $
   \mb{e}_{i,Au} = \mb{e}_i - \sum_{\mb{u}_j = \mb{\underline{u}}_j}\mb{b}_{ij}\mb{\underline{u}}_j - \sum_{\mb{u}_j = \mb{\bar{u}}_j} \mb{b}_{ij}\mb{\bar{u}}_j,$
    and $\alpha^\dagger_j = 0, j \in N\setminus\{i\}$.  If $\alpha_i^\dagger \geq 0$, there exists the unique solution to VI in $F(\mb{u}|\alpha^\dagger)$, denoted by $\mb{u}^\dagger$, which satisfies $\mb{b}_i^T \mb{u^\dagger} + \mb{e}_i - \bar{\mb{x}}_i = 0$.
    
\label{Theorem2}
\end{theorem}

    

\begin{proof}
    Based on the vector in \eqref{gradient}, we have that
    \begin{align}
    F(\mb{u} | \alpha^\dagger) = \mb{Q} \mb{u} + \alpha_i^\dag e^{\beta_i (\mb{x}_i - \mb{\bar{x}}_i)} \hat{\mb{b}}_i.
    \label{operator}
    \end{align}
    We compute the Jacobian of $F(\mb{u} | \alpha^\dagger)$ as
    \begin{align}
        J_F(\mb{u} |\alpha^\dagger) = \mb{Q} + \alpha_i^\dag \beta_i e^{\beta_i (\mb{x}_i - \mb{\bar{x}}_i)} \hat{\mb{b}}_i \mb{b}_i^T.
        \label{hessian}
    \end{align}
    Since matrix $ \hat{\mb{b}}_i \mb{b}_i^T$ is rank-1, we obtain the only non-zero eigenvalue as $\mb{b}_i^T\hat{\mb{b}}_i >0$. Combine with $\alpha_i^\dagger \geq0$ and Assumptions 1-2, we have that $J_F(\mb{u} |\alpha^\dagger)$ is positive-definite over $\mb{u} \in \mathbb{R}^{2n}$. According to \textit{Proposition 12.3} in \cite{rockafellar1998variational}, $F(\mb{u} |\alpha^\dagger)$ is strictly monotone. Thus, there exists the unique solution such that the VI to $J_F(\mb{u} |\alpha^\dagger)$ holds (\textit{Proposition 12.54} in \cite{rockafellar1998variational}).
    \par The subsequent proof is to show that the solution of $\mb{u}_{Au}$ to \eqref{def1} with other saturated actions is the solution $\mb{u}^\dagger$.
    \par We first expand it and get $\sum_{j=1}^{2n}F(\mb{u} | \alpha^\dagger)_j (\mb{u} - \mb{u}^\dagger)_j \geq 0$,
    where $(\cdot)_j$ represents the $j$-th element of vector. It suffices to prove that each term is non-negative. Since $\alpha_i^\dagger \geq 0$, neither active nor reactive power can hit the upper boundaries. To see this, assume that $\mb{u}^\dagger_j = \bar{\mb{u}}_j$ for some $j$. We get $ F(\mb{u} | \alpha^\dagger)_j > 0$, which violates \eqref{converged point1}. Thus, we proceed the proof by considering the following two cases:
    (a) If $\mb{u}^\dagger_j\in (\mb{\underline{u}}_j, \mb{\bar{u}}_j)$, by substituting the second row of linear equation \eqref{def1} into \eqref{operator}, we require that $F(\mb{u}^\dagger | \alpha^\dagger)_j = (\mb{Q}\mb{u}^\dagger + \alpha^\dagger_i \hat{\mb{b}}_i)_j = 0$, which is satisfied by the first part of \eqref{def1}; (b) If $\mb{u}_j^\dagger=\mb{\underline{u}}_j$, we require that 
    $F(\mb{u}^\dagger | \alpha^\dagger)_j = (\mb{Q}\mb{u}^\dagger + \alpha^\dagger_i \hat{\mb{b}}_i)_j \geq 0$.
    Consider the unsaturated case of \eqref{def1}, i.e. $A_u = 2N$, we have $ F(\tilde{\mb{u}}^\dagger | \tilde{\alpha}^\dagger)_j = 0, \: \forall j \in 2N$. Since $\hat{\mb{b}}_i$ is element-wise positive and $\mb{Q}$ is diagonal, the cut-off value caused by lower bound is shifted to unsaturated variables. Thus, we have that $\alpha_i^\dagger \geq \tilde{\alpha}_i^\dagger$. Combine with $\tilde{\mb{u}}_j^\dagger \leq \mb{\underline{u}}_j = \mb{u}^\dagger_j$, we obtain that $$(\mb{Q}\mb{u}^\dagger + \alpha^\dagger_i \hat{\mb{b}}_i)_j \geq (\mb{Q}\tilde{\mb{u}}^\dagger + \tilde{\alpha}_i^\dagger \hat{\mb{b}}_i)_j = 0.$$
    \par Above all, we prove that the solution of $\mb{u}_{Au}$ to \eqref{def1} with other saturated actions is the unique solution to VI to $F(\mb{u} | \alpha^\dagger)$ Then, we conclude that $\mb{b}_i^T \mb{u^\dagger} + \mb{e}_i - \bar{\mb{x}}_i = 0$.
\end{proof}

\par We further study the safety guarantee at the convergence of Algorithm \ref{alg} with barrier parameters $\alpha^s \in \mathbb{R}^n$, denoted by $\mb{u}^s$. We set $\alpha^s$ as $\alpha^s_i = \hat{\alpha}_i + \gamma^s$ at  $i = \arg\max_j \mb{x}_j$ with some $\gamma^s \geq 0$, and $\alpha^s_j = 0, j \in N\setminus\{i\}$, where $\hat{\alpha}_i$ is the solution to 
    \begin{align}
    \begin{bmatrix}
    \mb{\mb{Q}}_{A_u} & \mb{\hat{b}}_{{i,A_u}} \\
    \mb{\hat{b}}_{i,A_u}^T & \mb{0}
    \end{bmatrix}
    \begin{bmatrix}
    \mb{u}_{A_u} \\
    \alpha_{i}
    \end{bmatrix}
    =
    \begin{bmatrix}
    \mb{0}\\
    \mb{\bar{x}}_{i} - \mb{\hat{e}}_{i,A_u}
    \end{bmatrix}.
    \label{def2}
    \end{align}
\par In addition, similar to \eqref{converged point1}, we have $\mb{u}^s$ satisfy 
\begin{align}
    \langle F(\mb{u}^s | {\alpha}^s), \mb{u} - \mb{u}^s \rangle \ge 0, \quad \forall \mb{u} \in [\underline{\mb{u}}, \bar{\mb{u}}].
    \label{converged point2}
\end{align}
Since we focus on overvoltage problem, our safety guarantee in Theorem \ref{Theorem3} is built on $\mb{\bar{u}} = 0$.

\begin{theorem}
(Safety guarantee) Under Assumptions 1-2 and inaccurate model $\mb{\hat{B}}$, if upper bounds of actions satisfy $\bar{\mb{u}} = 0$, when \textbf{Algorithm \ref{alg}} under barrier weights $\alpha^s \in \mathbb{R}^n$ and fixed $\mb{\beta} \in \mathbb{R}_+^{n}$ achieves the convergence, denoted by $\mb{u}^s$, we have $\alpha_i^s \geq 0$, at $i=\arg\max_j \mb{b}_j^T \mb{u^s} + \mb{e}_j$ and 
\begin{align}
    \mb{b}_i^T \mb{u^s} + \mb{e}_i - \bar{\mb{x}}_i \leq 0,
    \label{safety}
\end{align}
with $\mb{\|B - \hat{B}\|} \leq \epsilon_B$ and 
\begin{align}
\gamma^s = \frac{\epsilon_B}{|\mb{\hat{b}}_{i, Au}^T \mb{Q}_{Au}^{-1}\mb{\hat{b}}_{i,Au}|} (\|\mb{\underline{u}}\| + \|\mb{Q}_{Au}^{-1}\mb{\hat{b}}_{i, Au}\| |\hat{\alpha}_i|).
\label{gamma}
\end{align}
\label{Theorem3}
\end{theorem}

\begin{proof}
    We proceed the proof by discussing different attention nodes at convergence $\mb{u}^s$. First, if $\mb{e}_i <\mb{\bar{x}}_i$ at $i = \arg\max_j \mb{b}_j^T \mb{u}^s + \mb{e}_j$, with \eqref{def1}, the maximum voltage deviation is achieved at the node without initial violation. Since $\bar{\mb{u}} = 0$,  \eqref{safety}  holds.
    \par Then by contradiction, we assume that $\alpha_i^s <0$. At the convergence $\mb{u}^s$, we have $F(\mb{u}^s | \alpha^s) < 0$. Based on \eqref{converged point2}, $\mb{u}^s$ can only locate on the upper bounds $\bar{\mb{u}} = 0$. With the activation of algorithm \ref{alg}, there exists node that incurs the voltage violation, which contradicts $\eqref{safety}$. Then we conclude that $\alpha_i^s \geq 0$ at $i=\arg\max_j \mb{b}_j^T \mb{u^s} + \mb{e}_j$.
    \par If  $\mb{e}_i  \geq \mb{\bar{x}}_i$, as  $\bar{\mb{u}} = 0$, the solution of $\alpha_i^\dagger$ to \eqref{def1} is non-negative. Based on Theorem \ref{Theorem2}, there exists $\alpha^\dagger \in \mathbb{R}^n_+$ with corresponding unsaturated set $A_u$ and attention node $i$, such that the unique solution to VI in $F(\mb{u}| \alpha^\dagger)$, denoted by $\mb{u}^\dagger$, satisfies 
    $\mb{b}_i^T \mb{u}^\dagger + \mb{e}_i - \bar{\mb{x}}_i = 0$. 
    \par Then, we follow a standard first-order perturbation analysis on \eqref{def1} and \eqref{def2} to obtain the variation of barrier weights $\alpha_i^\dagger$ and $\hat{\alpha}_i$. First, we calculate the error of coefficient as follows,
    \begin{align}
    |(\mb{b}_{i,Au} - \mb{\hat{b}}_{i,Au})^T \mb{Q}_{Au}^{-1}\mb{\hat{b}}_{i,Au}| &\leq \|\mb{b}_{i,Au} - \mb{\hat{b}}_{i,Au}\| \|\mb{Q}_{Au}^{-1}\mb{\hat{b}}_{i,Au}\| \nonumber \\
    &\leq \epsilon_B  \|\mb{Q}_{Au}^{-1}\mb{\hat{b}}_{i,Au}\|.
    \end{align}
    Consider the voltage drop with saturation, we have $\mb{\hat{e}}_{i,Au} - \mb{e}_{i,Au} = \sum_{\mb{u}_j = \mb{\underline{u}}_j} (\mb{b}_{ij} - \mb{\hat{b}}_{ij}) \mb{\underline{u}}_j.$ Then we obtain the error bound on constant as $|\mb{\hat{e}}_{i,Au} - \mb{e}_{i,Au}| \leq \epsilon_B \|\mb{\underline{u}}\|$.
    Combining the derived bounds, we obtain $|\hat{\alpha}_i - \hat{\alpha}_i^\dag| \leq \gamma^s$, with $\gamma^s$ given in \eqref{gamma}. With $\alpha_i^s = \hat{\alpha}_i + \gamma^s$, we have that $\alpha_i^s \geq \alpha_i^\dag \geq 0$.
    \par We proceed the proof by contradiction. At $i = \arg\max_j \mb{b}_j^T \mb{u}^s + \mb{e}_j$, assume that 
    \begin{align}
        \mb{b}_i^T \mb{u}^s + \mb{e}_i - \bar{\mb{x}}_i> 0.
        \label{assumption}
    \end{align}     
     Consider \eqref{def1} and \eqref{def2}, if $\mb{u}_j^s \in (\mb{\underline{u}}_j, \bar{\mb{u}}_j)$, we have  $F(\mb{u}^s | \alpha^s)_j = F(\mb{u}^\dagger | \alpha^\dagger)_j = 0$, i.e., 
    \begin{align*}
    &(\mb{Q} \mb{u}^s + \alpha_i^s e^{\beta_i(\mb{b}_i^T\mb{u}^s + \mb{e}_i - \mb{\bar{x}}_i)} \hat{\mb{b}}_i)_j \\
    &= (\mb{Q} \mb{u}^\dagger + \alpha_i^\dagger e^{\beta_i(\mb{b}_i^T\mb{u}^\dagger + \mb{e}_i - \mb{\bar{x}}_i)} \hat{\mb{b}}_i)_j = 0,
    \end{align*}
    With the assumption in \eqref{assumption} with $\hat{\mb{b}}_i \geq 0$ and diagonal $\mb{Q}$, we have $\mb{u}_j^s \leq \mb{u}_j^\dagger$. For $\mb{u}_j^s = \mb{\underline{u}}_j$, we have $\mb{u}_j^s = \mb{u}_j^\dagger$. Above all, we get that $\mb{b}_i^T \mb{u}^s + \mb{e}_i - \mb{\bar{x}}_i \leq \mb{b}_i^T \mb{u}^\dag + \mb{e}_i - \mb{\bar{x}}_i=0 $, which contradicts \eqref{assumption}. Thus, we conclude that $\mb{b}_i^T \mb{u^s} + \mb{e}_i - \bar{\mb{x}}_i \leq 0$, i.e., voltages at the convergence locate within the safety limit.
\end{proof}

At last, we show that Algorithm \ref{alg} achieves the exponential asymptotic convergence rate with proper choice of step-size.
\begin{theorem}
    (Exponential convergence)  Under the barrier parameter $\alpha^s \in \mathbb{R}^n$ in Theorem \ref{Theorem3} and upper bounds $\bar{\mb{u}} = 0$, the sequence $\{\mb{u}(k)\}$ generated by $\mb{u}(k+1) = \mb{u}(k) - \eta F(\mb{u}(k) | \alpha^s)$ has the exponential convergence rate near minimizer $\mb{u}^s$, i.e., 
    \begin{align}
   \frac{\|\mb{u}(k+1) - \mb{u}^s\|}{\|\mb{u}(k) - \mb{u}^s\|} \leq \rho,
   \label{conv0}
    \end{align}
    where $\rho = \max\{|1- \eta m_1|, |1- \eta m_2|\}$, $m_1 = \lambda_{\min}(\mb{Q}_{Au})$ and $m_2 = \lambda_{\max} (\mb{Q}_{Au}) + \alpha_i^s \beta_i \|\hat{\mb{b}_i}\|(\|\hat{\mb{b}_i}\| + \epsilon_B)$, with $\lambda_{\min}$ and $\lambda_{\max}$ as minimum and maximum eigenvalues, and $i = \arg\max_j \mb{b}_j^T \mb{u^s} - \mb{e}_j$. The step-size $\eta$ is chosen as
    \begin{align*}
    \eta \in (0, \frac{2}{L_s}), \quad L_s = \lambda_{\max} (\mb{Q}_{Au}) + \sigma\alpha_i^s \beta_i \|\hat{\mb{b}_i}\|(\|\hat{\mb{b}_i}\| + \epsilon_B),
    \end{align*}
    where $\sigma = \sup_k \max_{j} e^{\beta_j (\mb{x}_j(k) - \mb{\bar{x}}_j)}$ around the $\mb{u}^s$.
\label{theorem3}
\end{theorem}

\begin{proof}
    With the voltage feedback from grids as $\mb{x}_j(k) = \mb{b}_j^T\mb{u}(k) + \mb{e}_j, \: j \in N$, we expand $F(\mb{u}(k) | \alpha^s)$ at the stationary point $\mb{u}^s$ over unsaturated actions
    \begin{align}
        F(\mb{u}(k) | \alpha^s) =  F(\mb{u}^s | \alpha^s) + J_F(\mb{u}^s | \alpha^s) (\mb{u}(k) - \mb{u}^s).
    \end{align}
    \par By substituting into the sequence $\{\mb{u}(k)\}$, we have that
    \begin{align}
        \mb{u}(k+1) - \mb{u}^s = (\mathbb{I} - \eta J_F (\mb{u}^s  | \alpha^s)) (\mb{u}(k) - \mb{u}^s).
    \end{align}
    Thus, the convergence property is govern by the maximum eigenvalue of $\mathbb{I} - \eta J_F (\mb{u}^s  | \alpha^s)$, denoted by  $\rho$.
    \par With Theorem \ref{Theorem3}, we can bound $J_F (\mb{u}^s | \alpha^s)$ as below 
    \begin{align}
        \mb{Q}_{Au} + \alpha_i^s \beta_i \hat{\mb{b}_i}\mb{b}_i^T \succeq J_F(\mb{u}^s  | \alpha^s) \succeq \mb{Q}_{Au},
    \end{align}
    where the first inequality follows Theorem \ref{Theorem3}, i.e., $\max_j \mb{b}_j^T \mb{u^s} + \mb{e}_j - \bar{\mb{x}}_j \leq 0 $ and the second is supported by $\alpha_i^s \geq 0$. By Weyl's inequality, we have 
    \begin{subequations}
        \begin{align}
            \nonumber
            \lambda_{\max}( J_F (\mb{u}^s  | \alpha^s)) &\leq  \lambda_{\max} (\mb{Q}_{Au}) + \alpha_i^s \beta_i \mb{b}_i^T\hat{\mb{b}_i} \\  
            &\leq  \lambda_{\max} (\mb{Q}_{Au}) + \alpha_i^s \beta_i \|\hat{\mb{b}_i}\|(\|\hat{\mb{b}_i}\| + \epsilon_B) \\
            \lambda_{\min}(J_F(\mb{u}^s  | \alpha^s))  &\geq \lambda_{\min}(\mb{Q}_{Au}).
        \end{align}
    \end{subequations}


    By setting $m_1$ and $m_2$ accordingly, we can obtain the exponential convergence in \eqref{conv0} with $\rho = \max\{|1- \eta m_1|, |1- \eta m_2|\}$.

    \par The step-size $\eta$ is determined by the Lipschitz constant of $ F(\mb{u}(k) | \alpha^s)$  around $\mb{u}^s$, as shown below,
    \begin{align*}
    \sup_{k} \|J_F(\mb{u}(k)  | \alpha^s)\| &= \sup_{k} \lambda_{\max}(J_F(\mb{u}(k)  | \alpha^s)) \nonumber  \\
    &\leq \lambda_{\max} (\mb{Q}_{Au}) + \sigma \alpha_i^s \beta_i \|\hat{\mb{b}_i}\|(\|\hat{\mb{b}_i}\| + \epsilon_B) = L_s
    \end{align*}
    where $\sigma = \sup_k \max_{j} e^{\beta_j (\mb{x}_j(k) - \mb{\bar{x}}_j)}$ and choosing  $\eta \in (0, \frac{2}{L_s})$ suffices the convergence. 
       \end{proof}
 

\section{Simulation results}
\label{Sec:sim}

We evaluate the online exponential barrier method on a single-phase 56-bus radial network with line parameters provided by \cite{farivar2012optimal}. Under the inaccurate models, we demonstrate its effectiveness in the optimal and safe voltage control, by showing the tiny optimality gap of regulated voltages and little violations of safety limits during intermediate steps, respectively. 

\subsection{Environment Setup}
\par We assume that PV inverters are randomly installed at a subset of buses.
Inverters are enabled to reduce the active power output from maximum value $\mathbf{p_{av}}$ to $\mathbf{0}$, and either produce or absorb the reactive power to some limit, which we choose $40\%$ of the available solar. Then, we have the saturation bounds as $\mathbf{u^p} \in [-\mathbf{p_{av}}, \mathbf{0}]$ and $\mathbf{u^q} \in [-0.4 \mathbf{p_{av}}, 0.4\mathbf{p_{av}}]$. For buses without inverters, both the upper and lower bounds are set to zero.

\par The magnitude of nominal voltage at each bus is $12kV$, with the safety requirements as $\pm 5 \%$.
We jointly normalize the original loads in \cite{farivar2012optimal} and available solar energy to simulate scenarios with overvoltage problem. The cost function is set to
$c(\mathbf{u}) = \sum_{i\in N}c_p\|u^p_i\|^2 + c_q\|u^q_i\|^2$ with $c_p = 3$ and $c_q = 1$.

\par The imprecise model estimations $\mb{\hat{B}}_1$ and $\mb{\hat{B}}_2$ are randomly generated by considering both parametric and topological inaccuracies, i.e., measurement errors in line impedance and permutations of nodal positions.
And we obtain the relative model errors as $\epsilon_B / \|\mb{B}\|= 44.1\%$ and $52.8\%$, respectively.



\par We compare our framework with three control groups, (i) the case of no control, (ii) directly solving the LCQP in \eqref{2-3}-\eqref{2-5} with correct and erroneous model estimations, and (iii) regularized online primal-dual method proposed in \cite{dall2016optimal}.

\subsection{Results on optimality}

\begin{figure}[htbp]
  \centering
  \includegraphics[width=0.8\linewidth]{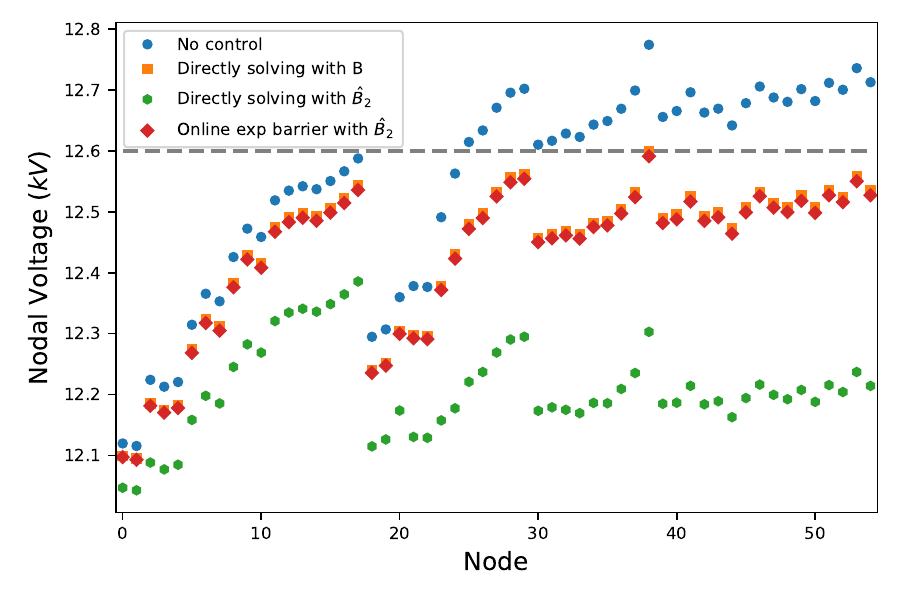}
    \caption{\footnotesize Comparison on the regulated voltages under online barrier method (red diamonds) and directly solving the LCQP (green hexagons) with inaccurate model $\mb{\hat{B}}_2$. The proposed method yields a solution close to that of the fully known model (orange squares), while no control (blue circles) exhibits significant deviations. The safety limit is denoted by the gray dashed line.}
\label{voltage2}
\end{figure}

Fig. \ref{voltage2} shows the capability of the proposed algorithm to control voltages within the safety limits under imprecise model estimations. In particular, the maximum nodal voltage is regulated to $12.59kV$ by implementing  Algorithm \ref{alg}, while the case of no control incurs the overvoltage at $12.77kV$ and directly solving the LCQP with $\hat{\mb{B}}_2$ results in the conservative actions at $12.39kV$. The corresponding barrier parameters are listed in Table \ref{table} and step-size is set as $\eta = 0.01$, calculated by Theorem \ref{Theorem3} and Theorem \ref{theorem3}, respectively.

\begin{table}[h]
	\renewcommand{\arraystretch}{1.3}
	\centering
	\begin{tabular}{>{\centering\arraybackslash}m{10em}|>{\centering\arraybackslash}m{6em}|>{\centering\arraybackslash}m{6em}}
		\ChangeRT{1.2pt}
		Inaccurate Models & $\mb{\hat{B}}_1$ &$\mb{\hat{B}}_2$ \\
        \ChangeRT{0.3pt}
		Safety Factor $\gamma^s$ &$2.7 \times 10^{-3}$ & $4.5 \times 10^{-3}$\\
        \ChangeRT{0.3pt}
		Barrier Weight $\alpha^s$ &$3.5 \times 10^{-3}$ &$5.4 \times 10^{-3}$\\
		\ChangeRT{1.2pt}
	\end{tabular}
	\caption{Initial safety factors and barrier weights of inaccurate model estimations $\mb{\hat{B}}_1$ and $\mb{\hat{B}}_2$, with $\beta_i = 200, \: i\in N$.}
	\label{table}
\end{table}

\subsection{Results on safety}

\begin{figure}[htbp]
  \centering

  \begin{subfigure}[b]{0.8\linewidth}
    \centering
    \includegraphics[width=\linewidth]{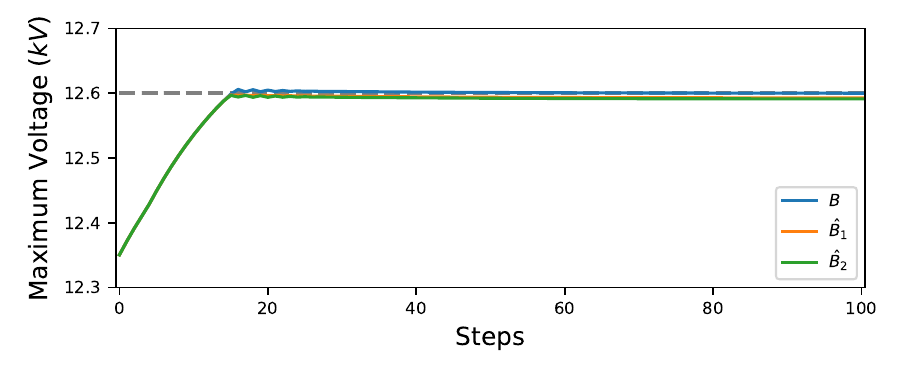}
    \caption{Online exponential barrier method}
  \end{subfigure}
  
  \vspace{0em}  

  \begin{subfigure}[b]{0.8\linewidth}
    \centering
    \includegraphics[width=\linewidth]{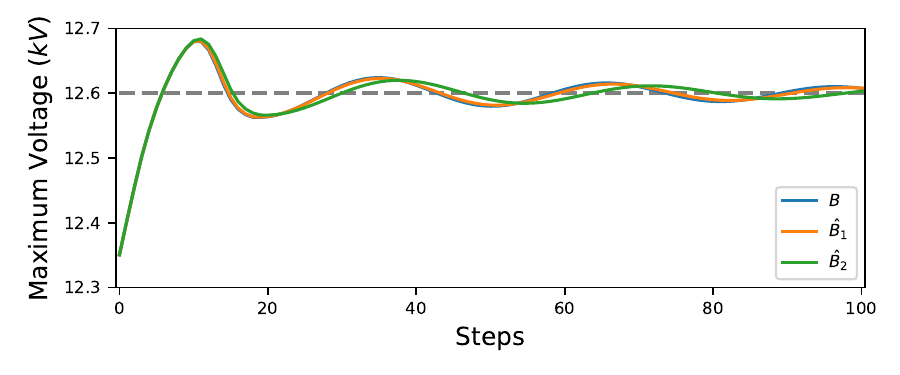}
    \caption{Online primal-dual method}
  \end{subfigure}

  \caption{\footnotesize Comparison on the maximum voltage during intermediate steps under (a) online exponential barrier method and (b) online primal-dual method, with models $\mb{B}$, $\mb{\hat{B}}_1$ and $\mb{\hat{B}}_2$. Although both are robust to the model inaccuracies, (b) exhibits significant violations in the process.}
  \label{voltage1}
\end{figure}

Fig. \ref{voltage1} exhibits the satisfaction of voltage safety limits during the interactions with physical grids. Although both frameworks are robust to model errors, online primal-dual method violates the constraints intermediately,
Moreover, starting from the initial action with $\kappa = 0.6$, our method converges in less than 20 steps, compared to more than 100 steps needed by the online primal-dual method.



\section{Conclusion}
\label{Sec:conclu}

In this paper, we proposed the online exponential barrier method for the optimal and safe voltage regulation in the presence of inaccurate distribution network models. Based on the integration of operational constraints as exponential barriers, we leveraged the online grids interactions to drive voltages to the proximity of optimum while achieving little violations of safety limits in the process. 
We provided the closed-form principles of barrier parameter selection and safety guarantees at the convergence. 
The performance on the optimality and safety was validated by simulations on the 56-bus radial test feeder. 
Future work includes the exploration of global convergence analysis and generalization to time-varying profiles of loads and solar energy with tracking properties. 


\section{Acknowledgment}
The authors would like to thank Manish. K. Singh for valuable discussions that helped improve this work.


\bibliographystyle{IEEEtran}
\bibliography{ref}

\end{document}